\documentclass[a4paper]{amsart} 
\usepackage[english]{babel} 
\usepackage[latin1]{inputenc}  
\usepackage[T1]{fontenc}       
\usepackage[paper=a4paper]{geometry}
\usepackage{url} 
\usepackage{amsmath,amssymb,latexsym,amsthm}
\usepackage{mathrsfs}
\usepackage{color} 
\usepackage{epsfig}
\usepackage{graphicx}
\usepackage{desclist}
\usepackage{enumerate, color}
\usepackage{hyperref}

\usepackage{pgf,tikz}
\usetikzlibrary{positioning,shapes,decorations.markings}
\usetikzlibrary{arrows}
\usetikzlibrary[patterns]

\input xy
\xyoption{all}

\newtheorem{theorem}{Theorem}[section]
\newtheorem{definition}[theorem]{Definition}

\newtheorem{proposition}[theorem]{Proposition}
\newtheorem{lemma}[theorem]{Lemma}
\newtheorem{remark}[theorem]{Remark}

\newcommand{\R}{\mathbb{R}}

\newcommand{\N}{\mathbb{N}}

\newcommand{\C}{\mathbb{C}}

\newcommand{\D}{\mathbb{D}}

\let\sm=\smallskip

\let\te=\textstyle

\let\phe=\varphi

\def\Hol{\mathop{\rm Hol}\nolimits}
\def\Fix{\mathop{\rm Fix}\nolimits}
\def\mlog{\mathop{{\te\frac{1}{2}\log}}}

\let\eps=\varepsilon \let\phe=\varphi
\let\de=\partial

\def\id{\mathop{\rm Id}\nolimits}

\def\1#1{\overline{#1}}
\def\2#1{\widetilde{#1}}
\def\3#1{\widehat{#1}}
\def\4#1{\mathbb{#1}}
\def\5#1{\frak{#1}}
\def\6#1{{\mathcal{#1}}}
\def\7#1{{\bf{#1}}}

\def\and{{\quad\text{and}\quad}}

\def\frak{\mathfrak }

\begin{document}
\title{Backward iteration in strongly convex domains.\\ \textit{Errata corrige}}
\author[Marco Abate]{Marco Abate}
\address{Marco Abate\\ Dipartimento di Matematica\\ Universit\`a di Pisa\\ Largo Pontecorvo 5, 56127 Pisa\\ Italy.} \email{marco.abate@unipi.it}
\author[Jasmin Raissy]{Jasmin Raissy*}
\address{Jasmin Raissy\\ Institut de Math\'ematiques de Toulouse, UMR5219\\ Universit\'e de Toulouse, CNRS\\ UPS IMT, 118 route de Narbonne, F-31062 Toulouse Cedex 9\\ France} 
\email{jraissy@math.univ-toulouse.fr}
%

\begin{abstract}
We correct a gap in two lemmas in \cite{AbaRai}, providing a new proof of the main results of that paper.
\end{abstract}

\maketitle

We have found a gap in the proofs of Lemmas 2.2 and 2.5 of our paper \cite{AbaRai}. In this note we fill these gaps, giving a proof of the main results using different arguments. For the benefit of the reader, we shall report here the complete proof of the main theorem, including the needed background results.

\section{Preliminary results}

The main tool we are going to use is the Kobayashi distance (see, e.g., \cite{Aba1} or \cite{Kob} for its main properties). In particular, we shall rely on the notion of complex geodesics, that we now recall. 

\begin{definition}
A {\em complex geodesic} in a hyperbolic manifold~$X$ is a holomorphic map
$\phe\colon\D\to X$ which is an isometry with respect to the Kobayashi distance of~$\D$
and the Kobayashi distance of~$X$, where $\D$ is the unit disk in the complex plane. 
\end{definition}

The following statements summarize the main results obtained by Lempert \cite{L} and Royden-Wong \cite{RW} on complex geodesics in strongly convex domains:


\begin{theorem}[{\cite[Theorem~2.6.19 and Corollary~2.6.30]{Aba1}}]
Let $D\Subset\C^n$ be a bounded convex domain. Then for every pair of distinct points $z$, $w\in D$ there exists a complex geodesic
$\phe\colon\D\to D$ such that $\phe(0)=z$ and $\phe(r)=w$, where $0<r<1$ is such that
$k_\D(0,r)=k_D(z,w)$; furthermore, if $D$ is strongly convex then $\phe$ is unique. Moreover a holomorphic map $\phe\in\Hol(\D,D)$ is a complex
geodesic if and only if $k_D\bigl(\phe(\zeta_1),\phe(\zeta_2)\bigr)=k_\D(\zeta_1,\zeta_2)$
for a pair of distinct points $\zeta_1$, $\zeta_2\in\D$.
\end{theorem}

\begin{proposition}[{\cite[Proposition~2.6.22]{Aba1}}] 
Let $D\Subset\C^n$ be a bounded convex domain. Then every complex geodesic $\phe\in\Hol(\D,D)$ admits
a {\em left inverse,} that is a holomorphic map $\widetilde p_\phe\colon D\to\D$ such that
$\widetilde p_\phe\circ\phe=\id_\D$. The map $p_\phe=\phe\circ\widetilde p_\phe$ is then 
a holomorphic retraction of $D$ onto the image of~$\phe$.
\end{proposition}

\begin{theorem}[{\cite[Theorem~2.6.29]{Aba1}}]
Let $D\Subset\C^n$ be a bounded strongly convex with~$C^2$ boundary. Then every
complex geodesics  $\phe$ extend continuously (actually, $\frac{1}{2}$-H\"older) to the boundary of~$\D$, and the image of $\phe$ is transversal to~$\de D$.
\end{theorem}

\begin{theorem}[{\cite[Theorem~2.6.45]{Aba1}}] 
Let $D\Subset\C^n$ be a bounded strongly convex with~$C^2$ boundary. Then 
for every
$z\in D$ and $\tau\in\de D$ there is a complex geodesic $\phe\in\Hol(\D,D)$ with
$\phe(0)=z$ and $\phe(1)=\tau$. Moreover for every pair of distinct points $\sigma$, $\tau\in\de D$
there is a complex geodesic $\phe\in\Hol(\D,D)$ such that $\phe(-1)=\sigma$ and
$\phe(1)=\tau$. 
\end{theorem}

The statement of \cite[Theorem~2.6.45]{Aba1} requires $D$ with~$C^3$ boundary, but the
proof of the existence actually works assuming just $C^2$ smoothness.

Now let $D\Subset\C^n$ be a bounded strongly convex domain with $C^2$ boundary, and $f\in\Hol(D,D)$ a holomorphic self-map of $D$. If the set $\Fix(f)$ of fixed
points of $f$ in~$D$ is not empty, then the sequence $\{f^{\circ k}\}$ of iterates of~$f$ is relatively
compact in~$\Hol(D,D)$, and 
there exists a submanifold $D_0\subseteq D$, the {\em limit manifold}
of~$f$, such that every limit of a subsequence of iterates is of the form $\gamma\circ\rho$,
where $\rho\colon D\to D_0$ is a holomorphic retraction, and $\gamma$ is a biholomorphism
of~$D_0$; furthermore, $f|_{D_0}$ is a biholomorphism of~$D_0$, and $\Fix(f)\subseteq D_0$ (see \cite{Aba01} or \cite[Theorem 2.1.29]{Aba1}). 

\begin{definition}
Let $D\Subset\C^n$ be a bounded strongly convex domain with $C^2$ boundary. A holomorphic map 
$f\in\Hol(D,D)$ is 
\begin{itemize}
\item {\em elliptic} if $\Fix(f)\ne\emptyset$,
\item {\em strongly elliptic} if its limit manifold
reduces to a single point, the {\em Wolff point} of the strongly elliptic map. 
\end{itemize}
We say that a point $p\in\Fix(f)$ is {\em attracting} if all the 
eigenvalues of $df_p$ have modulus strictly less than~1.
\end{definition}

We have an equivalent characterization of strongly elliptic maps:

\begin{lemma}[{\cite[Lemma 1.1]{AbaRai}}]
\label{uunoa}
Let $D\Subset\C^n$ be a bounded strongly convex domain with $C^2$ boundary, and $f\in\Hol(D,D)$.
Then the following facts are equivalent:
\begin{enumerate}[{\rm(i)}]
\item $f$ is strongly elliptic;
\item the sequence of iterates of $f$ converges to a point $p\in D$;
\item $f$ has an attracting fixed point $p\in D$;
\item there exists $p\in\Fix(f)$ such that $k_D\bigl(p,f(z)\bigr)< k_D(p,z)$ for all $z\in D\setminus\{p\}$.
\end{enumerate}
\end{lemma}

%
%

In the study of the dynamics of self-maps without fixed points, a crucial r\^ole is played by the
horospheres, a generalization of the classical notion of horocycle. 

It is a non-trivial fact (see, e.g., \cite[Theorem~2.6.47]{Aba1} or \cite{BP}) that for a 
bounded strongly convex domain with $C^2$ boundary $D\Subset\C^n$ the limit
$$
\lim_{w\to\tau}[k_D(z,w)-k_D(p,w)]
$$
exists for every $\tau\in\de D$ and $p\in D$ and we can therefore define $h_{\tau,p}\colon D\to\R^+$ as
$$
\mlog h_{\tau,p}(z)=\lim_{w\to\tau}[k_D(z,w)-k_D(p,w)]\;.
$$
Then we can introduce the following definition:

\begin{definition}
Let $D\Subset\C^n$ be a 
bounded strongly convex domain with $C^2$ boundary. The {\em horosphere} of {\em center} $\tau\in\de D$,
{\em radius} $R>0$ and {\em pole} $p\in D$ is the set
$$
E_p(\tau,R)=\{z\in D\mid h_{\tau,p}(z)< R\}\;.
$$
\end{definition}

We shall need the following fact:

\begin{lemma}[{\cite[Lemma~2.7.16]{Aba1}}]\label{th:prho}
Let $D\Subset\C^n$ be a bounded strongly convex domain with $C^2$ boundary, and $\phe\in\Hol(\D,D)$ a complex geodesic. Put $p=\phe(0)$ and $\tau=\phe(1)$. Then
$$
\tilde p_\phe\bigl(E_p(\tau,R)\bigr)=E^\D_0(1,R)
$$
for any $R>0$, where $E^\D_0(1,R)$ is the horosphere of center~$1$, pole~$0$ and radius~$R$ in~$\D$.
\end{lemma}

We can also introduce {\em $K$-regions} in a similar way. 

\begin{definition}\label{K-region}
Let $D\Subset\C^n$ be a 
bounded strongly convex domain with $C^2$ boundary.  The {\em $K$-region $K_p(\tau,M)$ 
of center $\tau\in\de D$, amplitude $M>0$ and pole $p\in D$} is the set
$$
K_p(\tau,M)=\{z\in D\mid \mlog h_{\tau,p}(z)+k_D(p,z)< \log M\}\;.
$$
\end{definition}

It is well-known (see \cite{Aba1,Ab2}) that the $K$-regions
with pole at the origin in the unit disk coincide with the classical Stolz regions,
and that the $K$-regions with pole at the origin in the unit ball $\mathbb{B}^n\subset\C^n$ coincide
with the usual Kor\'anyi approach regions. 

\begin{remark}\label{rem:admreg}
In strongly convex domains
$K$-regions are comparable to {\em Stein admissible approach regions $A(\sigma,M)$ of vertex $\sigma \in\de D$ and aperture $M>1$}:
\begin{equation}\label{admissiregion}
A(\sigma, M) 
=
\{z\in D \mid \|z-\sigma\|^2< M d(z,\de D), |\langle z-\sigma,n_\sigma\rangle| < M d(z,\de D)\}\;,
\end{equation}
where $n_\sigma$ is the outer unit normal vector to $\de D$ at $\tau$. Here ``comparable" means that for every $\tau\in\de D$ there exists a neighbourhood $U\subset\C^n$ of $\tau$ such that for any $M>1$ and $p\in D$ there are $M_1$, $M_2>1$ such that
$$
A(\sigma,M_1)\cap U\subseteq K_p(\sigma,M)\cap U\subseteq A(\sigma,M_2)\cap U\;;
$$
see, e.g., \cite[Propositions 2.7.4, 2.7.6 and p. 380]{Aba1}.
Moreover, changing the pole does not
change much the $K$-regions, because for each $p$,~$q\in D$ there is $L>0$ such that
\begin{equation}\label{eqKpole}
K_p(\tau, M/L)\subseteq K_q(\tau,M)\subseteq K_p(\tau,ML)
\end{equation}
for every $M>0$ (see \cite[Lemma 2.7.2]{Aba1}).
\end{remark}

\begin{definition}
Let $D\Subset\C^n$ be a 
bounded strongly convex domain with $C^2$ boundary.
Given $\tau\in\de D$, we shall say that a function $F\colon D\to\C^n$ has {\em $K$-limit
$\ell\in\C^n$ at~$\tau$} if $F(z)\to\ell$ as $z\to\tau$ inside any $K$-region centered at~$\tau$.
\end{definition}

Notice that the choice of the pole is immaterial because of \eqref{eqKpole}.

\begin{remark}\label{rem:Klimalim}
Thanks to Remark~\ref{rem:admreg}, a function $F\colon D\to\C^n$ has $K$-limit
$\ell\in\C^n$ at~$\tau$ if $F(z)\to\ell$ as $z\to\tau$ inside any Stein admissible region with vertex~$\tau$. In particular, the existence of a $K$-limit always
implies the existence of a non-tangential limit; see also \cite[Lemma~2.7.12]{Aba1} for a direct proof.
\end{remark}

Finally, the intersection of a horosphere (or $K$-region) of center~$\tau\in\de D$ and pole~$p\in D$ with the image of a complex geodesic~$\phe$
with $\phe(0)=p$ and $\phe(1)=\tau$ is the image via~$\phe$ of the horosphere (or $K$-region)
of center~1 and pole~0 in the unit disk (\cite[Proposition~2.7.8 and Lemma~2.7.16]{Aba1}).

\smallskip

We shall also need a generalization of the one-variable notion of angular derivative given by the
{\em dilation coefficient} (see \cite[Section~1.2.1 and Theorem~2.7.14]{Aba1})

\begin{definition}\label{BFP}
Let $D\Subset\C^n$ be 
a bounded strongly convex domain with $C^2$ boundary, $f\in\Hol(D,D)$, and let $\sigma\in\de D$. The {\em dilation coefficient $\beta_{\sigma,p}\in(0,+\infty]$
of $f$ at~$\sigma\in\de D$ with pole~$p\in D$} is given by
$$
\mlog\beta_{\sigma,p}=\liminf_{z\to\sigma}\bigl[k_D(p,z)-k_D\bigl(p,f(z)\bigr)\bigr]\;.
$$
Furthermore, $\sigma\in\de D$ is called a {\em boundary fixed point}
of~$f$ if $f$ has $K$-limit $\sigma$ at~$\sigma$. 
\end{definition}

Since
$$
k_D(p,z)-k_D\bigl(p,f(z)\bigr)\ge k_D\bigl(f(p),f(z)\bigr)-k_D\bigl(p,f(z)\bigr)\ge -k_D\bigl(p,f(p)\bigr)\;,
$$
the dilation coefficient cannot be zero. We also recall
the following useful formulas for computing the dilation coefficient obtained in \cite[Lemma~2.7.22]{Aba1}:
\begin{equation}\label{equsf}
\begin{aligned}
\mlog\beta_{\sigma,p}
&=\lim_{t\to 1}\bigl[k_D\bigl(p,\phe(t)\bigr)-k_D\bigl(p,f\bigl(\phe(t)\bigr)\bigr)
\bigr]\\
&=\lim_{t\to 1}\bigl[k_D\bigl(p,\phe(t)\bigr)-k_D\bigl(p,p_\phe\circ f\bigl(\phe(t)\bigr)\bigr)
\bigr]\;,
\end{aligned}
\end{equation}
where $\phe\in\Hol(\D,D)$ is a complex geodesic with $\phe(0)=p$ and $\phe(1)=\sigma$,
and $p_\phe=\phe\circ\widetilde p_\phe$ is the holomorphic retraction associated to~$\phe$.

When $\sigma$ is a boundary fixed point then the dilation coefficient does not depend
on the pole (see for example \cite[Lemma 1.3]{AbaRai}) and we shall then simply denote by $\beta_\sigma$ the dilation coefficient at a boundary fixed point.

\begin{definition}\label{BRFP}
Let $\sigma\in\de D$ be a boundary fixed point for a self-map $f\in\Hol(D,D)$ of
a bounded strongly convex domain with $C^2$ boundary $D\Subset\C^n$. We shall say that 
$\sigma$ is
\begin{itemize}
\item {\em attracting} if $0<\beta_\sigma<1$,
\item {\em parabolic} if $\beta_\sigma=1$,
\item {\em repelling} if $\beta_\sigma>1$.
\end{itemize}
\end{definition}

We can now quote the general version of Julia's lemma proved by Abate (see \cite[Proposition~2.4.15, Theorem~2.4.16 and Proposition~2.7.15]{Aba1}) that we shall need in this chapter.

\begin{proposition}[Abate, \cite{Aba1}]\label{Julia}
Let $D\Subset\C^n$ be a bounded strongly convex domain with $C^2$ boundary, and $f\in\Hol(D,D)$. Let $\sigma\in\de D$ and $p\in D$ be such that the dilation coefficient $\beta_{\sigma, p}$ is finite.
Then there exists a unique $\tau\in\de D$ such that for all $R>0$ we have
$$
f\bigl(E_p(\sigma,R)\bigr)\subseteq E_p\bigl(\tau,\beta_{\sigma,p}R\bigr)\;,
$$
and $f$ has $K$-limit $\tau$ at~$\sigma$. Moreover, if there is a sequence $\{w_\nu\}\subset D$ converging to $\sigma\in\de D$ so that $\{f(w_\nu)\}$ converges to $\tau_1\in\de D$ then $\tau=\tau_1$.
\end{proposition}

Finally, we recall the several variable version of the Wolff-Denjoy theorem given in~\cite{Aba01}
(see also \cite[Theorems~2.4.19 and~2.4.23]{Aba1} and \cite{Aba2}).

\begin{theorem}[Abate, \cite{Aba01}]\label{Wolff}
Let $D\Subset\C^n$ be a bounded strongly convex domain with $C^2$ boundary, and $f\in\Hol(D,D)$
without fixed points. Then there exists a unique $\tau\in\de D$ such that the sequence of
iterates of $f$ converges to~$\tau$. 
\end{theorem}

\begin{definition}\label{Wpoint}
Let $D\Subset\C^n$ be a bounded strongly convex domain with $C^2$ boundary, and $f\in\Hol(D,D)$
without fixed points. The point $\tau\in\de D$ introduced in the previous theorem is the
{\em Wolff point} of~$f$. 
\end{definition}

The dilation coefficient can also be used to characterize the Wolff point of $f\in\Hol(D,D)$ without fixed points in~$D$ defined above. 

\begin{proposition}[{\cite[Proposition~1.6]{AbaRai}}]
\label{AcW} 
Let $D\Subset\C^n$ be a bounded strongly convex domain with $C^2$ boundary, and $f\in\Hol(D,D)$ without fixed points in~$D$. Then the
following assertions are equivalent for a point $\tau\in\de D$:
\begin{enumerate}[{\rm(i)}]
\item $\tau$ is a boundary fixed point with $0<\beta_\tau\le 1$;
\item $f\bigl(E_p(\tau,R)\bigr)\subseteq E_p(\tau,R)$ for all $R>0$ and any (and hence all) $p\in D$;
\item $\tau$ is the Wolff point of $f$.
\end{enumerate}
\end{proposition}

\begin{definition}
Let $D\Subset\C^n$ be a bounded strongly convex domain with $C^2$ boundary, and $f\in\Hol(D,D)$
without fixed points and with Wolff point $\tau\in\de D$. 
We shall say that $f$ is {\em hyperbolic} if $0<\beta_\tau<1$
and {\em parabolic} if $\beta_\tau=1$.
\end{definition}

A final definition is needed. 

\begin{definition}\label{boundedstep}
Let $X$ be a Kobayashi hyperbolic manifold. We say that a sequence $\{z_k\}\subset X$ has
{\em bounded Kobayashi step} if 
$$
a=\sup\limits_k k_X(z_{k+1},z_k)<+\infty\;.
$$ 
The number~$a$ is
the {\em Kobayashi step} of the sequence.
\end{definition}

Notice that for any $f\in\Hol(X,X)$ and any $z_0\in X$ the orbit $\{f^k(z_0)\}$ has bounded Kobayashi step $k_X\bigl(z_0,f(z_0)\bigr)$.


\section{Main results}


\begin{definition}\label{backorbit}
Let $f\colon X\to X$ be a self-map of a set $X$. A {\em backward orbit} (or {\em backward iteration sequence\/}) for~$f$ is a sequence
$\{x_k\}_{k\in\N}\subset X$ so that $f(x_{k+1})=x_k$ for all $k\in\N$.
\end{definition}

The aim of this note is to prove the following version of \cite[Theorem~0.1]{AbaRai}:

\begin{theorem}
\label{OJMz:}
Let $D\Subset\C^n$ be a bounded strongly convex domain with $C^2$ boundary.
Let $f\in\Hol(D,D)$ be either hyperbolic or strongly elliptic, with Wolff point $\tau\in\overline{D}$. 
Let $\{z_k\}\subset D$ be a backward orbit for~$f$ with bounded Kobayashi step.
Then:
\begin{enumerate}[{\rm(i)}]
\item the sequence $\{z_k\}$ converges to a 
boundary fixed point~$\sigma\in\de D$;
\item if $\sigma\ne \tau$ then $\sigma$ is repelling;
\item $\sigma\ne\tau$ if and only if 
$\{z_k\}$ goes to~$\sigma$ inside a $K$-region, that is, there exists $M>0$ so that $z_k\in K_{p}(\sigma,M)$ eventually, where $p$ is any point in $D$. 
\end{enumerate}
\end{theorem}

\begin{remark}
If $f$ is strongly elliptic then clearly $\sigma\ne\tau$. We conjecture that $\sigma\ne\tau$ in the hyperbolic case too.
\end{remark}

\begin{remark}
The following proof does not work in the parabolic case, considered in the original version of \cite[Theorem~0.1]{AbaRai}. Thus the behavior of backward orbits for a parabolic self-map is still not understood, even (as far as we know) in the unit ball of $\C^n$ (see \cite{O}).
\end{remark}

\begin{proof} 
We begin with a first general lemma, saying that if a backward orbit with bounded Kobayashi step converges to a boundary point, then this point necessarily is a boundary fixed point.

\begin{lemma}[{\cite[Lemma 2.3]{AbaRai}}]\label{duetre}
Let $D\Subset\C^n$ be a bounded strongly convex domain with $C^2$ boundary, and let
$f\in\Hol(D,D)$. Let $\{z_k\}\subset D$ be a backward orbit for~$f$ with bounded Kobayashi step~$a=\mlog\alpha$ converging to~$\sigma\in\de D$. Then $\sigma$ is a boundary fixed point of~$f$ and 
$\beta_\sigma\le\alpha$. 
\end{lemma}

\begin{proof} 
Fix $p\in D$. 
First of all we have 
\begin{equation}
\begin{aligned}
\mlog\beta_{\sigma,p}=\liminf_{w\to\sigma}\bigl[k_D(w,p)-k_D\bigl(f(w),p\bigr)\bigr]&\le \liminf_{k\to+\infty}
[k_D(z_{k+1},p)-k_D\bigl(z_k,p\bigr)]\cr
&\le \liminf_{k\to+\infty} k_D(z_{k+1},z_k)\cr
&\le a=\mlog\alpha\;.
\end{aligned}
\end{equation}
Since $z_k\to\sigma$ and $f(z_k)=z_{k-1}\to\sigma$ as $k\to+\infty$, Proposition~\ref{Julia} yields the assertion.
\end{proof}

The rest of the proof is divided into two cases according to whether $f$ is hyperbolic or strongly elliptic. 

\sm\noindent{\bf Hyperbolic case.} 
In this case, we first prove that any backward orbit has to accumulate to the boundary of the domain $D$. 

\begin{lemma}[{\cite[Lemma 2.1]{AbaRai}}]\label{dueuno}
Let $D\Subset\C^n$ be a bounded strongly convex domain with $C^2$ boundary. Let $\{z_k\}\subset D$ be a backward orbit for a hyperbolic or parabolic self-map $f\in\Hol(D,D)$. 
Then $z_k\to\de D$ as $k\to+\infty$. 
\end{lemma}

\begin{proof}
Assume, by contradiction, that the sequence does not converge to $\de D$. Then there exists a subsequence $\{z_{k_n}\}$ converging to $w_0\in D$; in particular,  
$$
k_D(w_0,z_{k_n}) \to 0 \quad \hbox{as}~k_n\to +\infty\;.
$$
Therefore
$$
k_D\bigl(f^{k_n}(w_0), f^{k_n}(z_{k_n})\bigr) \le k_D(w_0,z_{k_n}) \to 0 \quad \hbox{as}~k_n\to +\infty\;.
$$
But, on the other hand, $f^{k_n}(z_{k_n}) = z_0$ for all $k_n$; moreover,
$f^{k_n}(w_0) \to\tau$ as $k_n\to+\infty$, where $\tau\in\de D$ is the Wolff point of~$f$. So
$$
\lim_{k_n\to\infty}k_D\bigl(f^{k_n}(w_0), f^{k_n}(z_{k_n})\bigr) = +\infty\;,
$$
because $k_D$ is complete, giving us a contradiction.
\end{proof}

In order to prove the convergence of the whole sequence towards a boundary fixed point $\sigma \in\de D$ we first need the following estimate.

\begin{lemma}[{\cite[Lemma 2.6]{AbaRai}}]\label{duesei}
Let $D\Subset\C^n$ be a bounded strongly convex domain with $C^2$ boundary, and fix $p\in D$.
Let $f\in\Hol(D,D)$ be hyperbolic or parabolic with Wolff point $\tau\in\de D$ and dilation coefficient
$0<\beta_\tau\le 1$. 
Let $\{z_k\}\subset D$ be a backward orbit for~$f$. Then for every $k\in\N$ we have
$$
h_{\tau,p}(z_k)\ge \left(\frac{1}{\beta_\tau}\right)^k h_{\tau,p}(z_0)\;.
$$
\end{lemma}

\begin{proof}
Put $t_k=h_{\tau,p}(z_k)$.
By definition, 
$z_k\in\de E_{p}(\tau,t_k)$. By Proposition~\ref{Julia}, if $z_{k+1}\in
E_{p}(\tau,R)$ then $z_k\in E_{p}(\tau,\beta_\tau R)$. Since $z_k\notin  E_{p}(\tau,t_k)$, we
have that $z_{k+1}\notin E_{p}(\tau, \beta_\tau^{-1}t_k)$, that is
\begin{equation}
t_{k+1}\ge \frac{1}{\beta_\tau}\, t_k\;,
\label{equuno}
\end{equation}
and the assertion follows by induction.
\end{proof}

This estimate allows us to prove part (i) in the hyperbolic case.

\begin{lemma}[{\cite[Remark 2.1]{AbaRai}}]\label{nuovoduedue}
Let $D\Subset\C^n$ be a bounded strongly convex domain with $C^2$ boundary.
Let $f\in\Hol(D,D)$ be hyperbolic with Wolff point $\tau\in\de D$ and let $\{z_k\}\subset D$ be a backward orbit for~$f$ with bounded Kobayashi step $a>0$. Then $\{z_k\}$ converges to a boundary fixed point $\sigma\in\de D$.
\end{lemma} 

\begin{proof}
First of all, recall that \cite[Lemma~2.4 and Remark~3]{AS} 
yields a constant $C_1>0$ such that
\begin{equation}
\|z_k-z_{k+1}\|^2 + |\langle z_k - z_{k+1}, z_k\rangle|\le \frac{C_1^2}{{1-\hat a^2}}{d(z_k,\de D)},
\label{eqASStima}
\end{equation}
and so
\begin{equation}
\|z_k-z_{k+1}\| \le \frac{C_1}{\sqrt{1-\hat a^2}}\sqrt{d(z_k,\de D)}\le \frac{C_1}{ 1-\hat a}\sqrt{d(z_k,\de D)}\;,
\label{eqAS}
\end{equation}
where $\hat a=\tanh a\in(0,1)$.
On the other hand, given $p\in D$
the triangular inequality and the upper estimate \cite[Theorem~2.3.51]{Aba1} on the boundary behaviour of the Kobayashi distance yield a constant $C_2>0$ such that
$$
\mlog h_{\tau,p}(z_k)\le k_D(p,z_k)\le C_2-\mlog d(z_k,\de D)\;,
$$ 
that is
\begin{equation}
d(z_k,\de D)\le \frac{e^{2C_2}}{ h_{\tau,p}(z_k)}\;,
\label{eqestup}
\end{equation}
and thus
\begin{equation}
\|z_k-z_{k+1}\|\le \frac{C}{ 1-\hat a}\sqrt{\frac{1}{ h_{\tau,p}(z_k)}}\;,
\label{equotto}
\end{equation}
for a suitable $C>0$. Therefore using  \eqref{equuno} we obtain that for every
$k,m\ge 0$ we have
\begin{equation}
\begin{aligned}
\|z_k - z_{k+m}\|
&\le\sum_{j=k}^{k+m-1}\|z_j - z_{j+1}\|
\le\frac{C}{ 1-\hat a}\frac{1}{\sqrt{h_{\tau,p}(z_k)}}\sum_{j=0}^{m-1} \beta_\tau^{j/2}\\
&\le\frac{C}{ 1-\hat a}\frac{1}{ 1-\beta_\tau^{1/2}}\frac{1}{\sqrt{h_{\tau,p}(z_k)}}\;.
\end{aligned}
\label{eqottobis}
\end{equation}
Since $h_{p.\tau}(z_k)\to+\infty$ as $k\to+\infty$ by Lemma~\ref{duesei} it follows that
$\{z_k\}$ is a Cauchy sequence in~$\C^n$, converging to a point $\sigma$, necessarily belonging to~$\de D$ by Lemma~\ref{dueuno}. The proof is then completed by quoting Lemma~\ref{duetre}.
\end{proof}

The following lemma allows us to control the dilation coefficient at the limit of a backward orbit, giving in particular part (ii) of Theorem~\ref{OJMz:} in the hyperbolic case.
 
\begin{lemma}[{\cite[Lemma 2.4]{AbaRai}}]\label{duequattro}
Let $D\Subset\C^n$ be a bounded strongly convex domain with $C^2$ boundary.
Let $f\in\Hol(D,D)$ be hyperbolic or parabolic with Wolff point $\tau\in\de D$
and dilation coefficient~$0<\beta_\tau\le 1$. Let $\sigma\in\de D\setminus\{\tau\}$ be a boundary fixed point with finite dilation coefficient~$\beta_\sigma$. 
Then
$$
\beta_\sigma\ge \frac{1}{ \beta_\tau}\ge1\;.
$$
In particular, if $f$ is hyperbolic then $\sigma$ is repelling.
\end{lemma}

\begin{proof} 
Let $\phe\colon\overline{\D}\to\overline{D}$ be a complex geodesic such that
$\phe(-1)=\sigma$ and $\phe(1)=\tau$, and set $p=\phe(0)$. 
Proposition~\ref{Julia} yields
$$
p\in \overline{E_p(\sigma,1)}\qquad\Longrightarrow\qquad f(p)\in\overline{E_p(\sigma,\beta_\sigma)}
$$
and 
$$
p\in \overline{E_p(\tau,1)}\qquad\Longrightarrow\qquad f(p)\in\overline{E_p(\tau,\beta_\tau)}\;.
$$
Hence $\overline{E_p(\sigma,\beta_\sigma)}\cap\overline{E_p(\tau,\beta_\tau)}\ne\emptyset$. 

Let $\widetilde p_\phe\colon D\to\D$ be the left-inverse of~$\phe$. Then using Lemma~\ref{th:prho} we get
$$
\emptyset
\ne \widetilde p_\phe\left(\overline{E_p(\sigma,\beta_\sigma)}\cap\overline{E_p(\tau,\beta_\tau)}\right)
	\!\!\subseteq\widetilde p_\phe(\overline{E_p(\sigma,\beta_\sigma)})\cap \widetilde p_\phe(\overline{E_p(\tau,\beta_\tau)})
	= \overline{E_0^\D(-1,\beta_\sigma)}\cap \overline{E_0^\D(1,\beta_\tau)}\;.
$$
Now, $E_0^\D(1,\beta_\tau)$ is an Euclidean disk of radius~$\beta_\tau/(\beta_\tau+1)$ tangent to~$\de\D$
in~$1$, and $E_0^\D(-1,\beta_\sigma)$ is an Euclidean disk of radius~$\beta_\sigma/(\beta_\sigma+1)$ tangent to~$\de\D$
in~$-1$. So these disks intersect if and only if
$$
1-\frac{2\beta_\tau}{ \beta_\tau+1}\le -1+\frac{2\beta_\sigma}{\beta_\sigma+1}\;,
$$
which is equivalent to $\beta_\sigma \beta_\tau\ge 1$, as claimed.
\end{proof}

We can now prove the first half of Theorem~\ref{OJMz:}.(iii) for the hyperbolic case.

\begin{lemma}
\label{tre}
Let $D\Subset\C^n$ be a bounded strongly convex domain with $C^2$ boundary. Let $f\in\Hol(D,D)$ be hyperbolic with Wolff point $\tau\in\de D$ and dilation coefficient~$0<\beta_\tau< 1$, and let $\{z_k\}\subset D$ be a
backward orbit with bounded Kobayashi step $a=\mlog\alpha$ converging to $\sigma\in\de D\setminus\{\tau\}$. Then for every $p\in D$ there exists $M>0$ such that $z_k\in K_{p}(\sigma,M)$ eventually.
\end{lemma}

\begin{proof}
Fix $p\in D$. By Remark~\ref{rem:admreg} it suffices to prove that there exists $M>1$ such that $\{z_k\}$ converges to~$\sigma$ inside an admissible approach region $A(\sigma,M)$.

Set again $t_k :=h_{\tau,p}(z_k)$. Thanks to \eqref{equuno}, we have
$$
\frac{1}{t_{k+m}} \le \beta_\tau^m \frac{1}{t_k}
$$
for all $k$, $m\ge 0$. Moreover, thanks to \cite[Corollary 2.3.55]{Aba1}, since $\sigma \ne \tau$, there exists $\varepsilon>0$ and $K>0$ such that for any $w\in D\cap B(\tau, \varepsilon)$ and $k\in\N$ such that $z_k\in D\cap B(\sigma, \varepsilon)$ we have 
$$
k_D(z_k,w)\ge - \mlog d(z_k,\de D) - \mlog  d(w,\de D) + K\;,
$$
where $B(x,\eps)$ is the Euclidean ball of center $x$ and radius~$\eps$.

On the other hand, \cite[Theorem 2.3.51]{Aba1} yields $c_1\in\R$ such that  
$$
k_D(w,p) \le c_1 -\mlog d(w,\de D)
$$
for any $w\in D$. So for $w\in D\cap B(\tau, \varepsilon)$ and $k$ sufficiently large we have 
$$
k_D(z_k,w)-k_D(w,p) \ge - \mlog d(z_k,\de D) - \mlog  d(w,\de D) +\mlog d(w,\de D) - c_1 + K\;,
$$
which implies 
$$
t_k=h_{\tau,p}(z_k) = \lim_{w\to\tau}[k_D(z_k,w)-k_D(w,p)] \ge - \mlog d(z_k,\de D) + K-c_1\;,
$$
that is 
\begin{equation}\label{7.12b}
\frac{1}{t_k} \le \widetilde C_1 d(z_k, \de D),
\end{equation}
for some $\widetilde C_1>0$.

Therefore, thanks to \eqref{eqottobis}, 
for all $m\ge 0$ and $k$ large enough we have 
\begin{equation}\label{stima1prima}
\|z_k - z_{k+m}\|
\le \frac{C \widetilde C_1}{ 1-\hat a}\frac{1}{ 1-\beta_\tau^{1/2}}  \sqrt{d(z_k, \de D)}
\end{equation}
for some $C>0$, and letting $m$ tend to infinity we obtain that for $k$ sufficiently large there is $M_1>1$ such that
\begin{equation}\label{stima1}
\|z_k-\sigma\|< M_1 \sqrt{d(z_k, \de D)}.
\end{equation}

On the other hand, up to translating the domain, without loss of generality we can assume that $D$ contains the origin. In particular, $D$ being bounded and strongly convex we can replace $n_\sigma$ by $\sigma$ in the definition of $A(\sigma,M)$. Therefore, to conclude the proof it suffices to prove that there exists $M_2>1$ such that 
$$
|\langle z_k-\sigma,\sigma\rangle| \le M_2 d(z_k, \de D)
$$
for $k$ large enough. Now
$$
|\langle z_j- z_{j+1}, z_j-\sigma \rangle| \le \|z_j- z_{j+1}\| \|z_j-\sigma\|,
$$
and so, thanks to \eqref{eqASStima} and \eqref{stima1}, for $k$ large enough and $m\ge 0$ we have
\begin{equation}\label{eqdieci}
\begin{aligned}
|\langle z_k - z_{k+m},\sigma\rangle|
&\le
\sum_{j=k}^{k+m-1} \left|\left\langle z_j - z_{j+1},\sigma\right\rangle\right|\\
&\le \sum_{j=k}^{k+m-1} \Big(|\langle z_j - z_{j+1}, z_j-\sigma \rangle| + |\langle z_j - z_{j+1}, z_j \rangle|\Big)\\
&\le
\sum_{j=k}^{k+m-1} \Big( \|z_j - z_{j+1}\| \|z_j-\sigma\| + \frac{C_1^2}{{1-\hat a^2}} d(z_j,\de D)\Big)\cr
&\le\sum_{j=k}^{k+m-1} \left( \frac{ M_1 C_1}{1-\hat a}{d(z_j,\de D)}+ \frac{C_1^2}{{1-\hat a^2}} d(z_j,\de D)\right)\\
&\le C' \sum_{j=k}^{k+m-1} d(z_j,\de D),
\end{aligned}
\end{equation}
for some $C'>0$. Arguing as in \eqref{eqottobis}, using \eqref{eqestup} and \eqref{7.12b} we obtain
$$
|\langle z_k - z_{k+m},\sigma\rangle|\le M_2 d(z_k,\de D)
$$
for $m\ge0$, $k$ large enough and for some $M_2>1$. Letting $m$ tend to infinity we finally have
$$
|\langle z_k - \sigma,\sigma\rangle|\le M_2 d(z_k,\de D).
$$
It then suffices to take $M= \max \{M_1, M_2\}$ to conclude the proof.
\end{proof}

The following lemma completes the proof of Theorem~\ref{OJMz:}.(iii):

\begin{lemma}
\label{lasttre}
Let $D\Subset\C^n$ be a bounded strongly convex domain with $C^2$ boundary. Let $f\in\Hol(D,D)$ be hyperbolic with Wolff point $\tau\in\de D$ and dilation coefficient~$0<\beta_\tau< 1$, and let $\{z_k\}\subset D$ be a
backward orbit with bounded Kobayashi step $a=\mlog\alpha$ converging to $\sigma\in\de D\setminus\{\tau\}$ inside a $K$-region. Then $\sigma\ne\tau$.
\end{lemma}

\begin{proof}
Assume, by contradiction, that $\sigma=\tau$. 
Fix $p\in D$, and let $M>1$ be such that $z_k\in K_p(\tau,M)$. Given $\eps>0$, \cite[Lemma 2.7.1]{Aba1} yields $r>0$ such that if $k_D(z_k,p)\ge r$ then 
$z_k\in E_p(\tau,\eps)$, that is $h_{\tau,p}(z_k)<\eps$. Since $k_D(z_k,p)\to+\infty$, it follows that $h_{\tau,p}(z_k)\to 0$ as $k\to+\infty$. But Lemma~\ref{duesei} implies that $h_{\tau,p}(z_k)\to+\infty$, contradiction.
\end{proof}

\sm\noindent{\bf Strongly elliptic case.} Also in this case, we start by proving by contradiction that any backward orbit has to accumulate to the boundary of the domain $D$. 

\begin{lemma}\label{duenove}
Let $D\Subset\C^n$ be a bounded strongly convex domain with $C^2$ boundary.
Let $f\in\Hol(D,D)$ be strongly elliptic with Wolff point $p\in D$, and let $\{z_k\}\subset D$ be a backward orbit with bounded Kobayashi step $a=\mlog\alpha$. Then $z_k\to\de D$ as $k\to+\infty$.
\end{lemma}

\begin{proof} Define $\ell_k>0$ by setting $\mlog \ell_k=k_D(z_k,p)$. Since $f$ is strongly elliptic, we have
$$
k_D(z_k,p)<k_D(z_{k+1},p)\;,
$$
and thus the sequence $\{\ell_k\}$ is strictly increasing. Assume, by contradiction, that it has a finite limit $\ell_\infty$. This means that every limit point $z_\infty$ of the sequence $\{z_k\}$ satisfies $k_D(z_\infty,p)=\mlog \ell_\infty$. But $f(z_\infty)$ is a limit point of the sequence $\{f(z_k)\}=\{z_{k-1}\}$ and thus we again have $k_D\bigl(f(z_\infty),p\bigr)=\mlog \ell_\infty$, which is impossible by Lemma~\ref{uunoa} because $f$ is strongly elliptic. Therefore $\ell_\infty=+\infty$, which means that $z_k\to\de D$.
\end{proof}

This allows us to prove the following key result.

\begin{lemma}\label{dueotto}
Let $D\Subset\C^n$ be a bounded strongly convex domain with $C^2$ boundary. Let $f\in\Hol(D,D)$ be strongly elliptic with Wolff point $p\in D$. Let $\{z_k\}\subset D$ be a backward orbit with bounded Kobayashi step. Then there exists a constant $0<c<1$ such that
$$
k_D(z_k,p)-k_D(z_{k+1},p)\le\mlog c<0
$$
for all $k\in\N$.
\end{lemma}

\begin{proof} Assume, by contradiction, that for every $0<c<1$ there is $k(c)\in\N$ such that 
$$
k_D(z_{k(c)},p)-k_D(z_{k(c)+1},p)>\mlog c\;,
$$
that is
$$
k_D(z_{k(c)+1},p)-k_D\bigl(f(z_{k(c)+1}), p\bigr)<-\mlog c\;.
$$
Consider the sequences $\{z_{k(1-\frac{1}{j})+1}\}$ and $\{z_{k(1-\frac{1}{ j})}=f(z_{k(1-\frac{1}{j})+1})\}$. Thanks to Lemma \ref{duenove}, we know that both these sequences accumulate on $\de D$; therefore, by extracting subsequences, we can find a subsequence $\{z_{k_j}\}$ such that $z_{k_j}\to \sigma_1\in \de D$, $f(z_{k_j})\to \sigma_2\in \de D$ as $j\to +\infty$ and
$$
\lim_{j\to+\infty}\bigl[k_D(z_{k_j},p)-k_D\bigl(f(z_{k_j}),p\bigr)\bigr]\le 0\;.
$$
If $\sigma_1\ne\sigma_2$, then \cite[Corollary 2.3.55]{Aba1}, together with the fact that $\{z_k\}$ has bounded Kobayashi step, lead to a contradiction since for $k$ large enough there is $K\in\R$ such that
$$
a
\ge
k_D\bigl(z_{k_j},f(z_{k_j})\bigr)\ge - \mlog d(z_{k_j},\de D) - \mlog  d\bigl(f(z_{k_ j}),\de D\bigr) + K
$$
whereas the right-hand side tends to infinity. Therefore, $\sigma_1=\sigma_2$ and we have
$$
\liminf_{z\to\sigma_1}\bigl[k_D(z,p)-k_D\bigl(f(z),p\bigr)\bigr]\le 0\;.
$$
Then we can apply Proposition~\ref{Julia}, obtaining that $\sigma_1$ is a boundary fixed point and that for any~${R>0}$ we have $f\bigl(E_p(\sigma_1,R)\bigr)\subseteq E_p(\sigma_1,R)$. 
We can then choose $R<1$ so that $p\notin\overline{E_p(\sigma_1,R)}$, and let $w\in\overline{E_p(\sigma_1,R)}$ be a point closest to~$p$ with respect to the Kobayashi distance. Since $f(w)\in\overline{E_p(\sigma_1,R)}$ this means that $k_D\bigl(f(w),p\bigr)\ge k_D(w,p)$, which is impossible because $w\ne p$ and $f$ is strongly elliptic.
\end{proof}

We can now prove that the whole backward orbit converges to a boundary fixed point $\sigma\in\de D$, which is obviously different from the Wolff point  $p\in D$.

\begin{lemma}[{\cite[Remark 2.2]{AbaRai}}]\label{duenovebis}
Let $D\Subset\C^n$ be a bounded strongly convex domain with $C^2$ boundary.
Let $f\in\Hol(D,D)$ be strongly elliptic with Wolff point $p\in D$, and let $\{z_k\}\subset D$ be a
backward orbit with bounded Kobayashi step $a=\mlog\alpha$. Then $\{z_k\}$ converges to a boundary fixed point $\sigma\in\de D$ with $\beta_\sigma\le\alpha$.
\end{lemma}

\begin{proof}
Without loss of generality, we can assume that $z_0\ne p$. We consider $s_k>0$ defined by setting $-\mlog s_k=k_D(z_k,p)$. Taking the constant $0<c<1$ given by the  Lemma~\ref{dueotto}, we therefore have
$$
-\mlog s_k+\mlog s_{k+1}\le\mlog c\;,
$$
that is
\begin{equation}\label{equsei}
s_{k+1}\le cs_k\;.
\end{equation}
Therefore $s_{k+m}\le c^m s_k$ for every $k,m\in\N$, and using again \eqref{eqASStima} and \cite[Theorem 2.3.51]{Aba1} as in the proof of Lemma \ref{nuovoduedue}, for all $j\in\N$ we obtain
$$
\|z_j-z_{j+1}\|\le \frac{C}{1-\hat a}\sqrt{s_j}
$$
for a suitable $C>0$, where $\hat a=\tanh a$. Arguing exactly as in \eqref{eqottobis} we then obtain
that 
\begin{equation}\label{eqnove}
\|z_k-z_{k+m}\| \le\frac{C}{1-\hat a}\frac{1}{1-c^{1/2}}\sqrt{s_k},
\end{equation}
for any $m\ge 0$ and $k$ large enough. So $\{z_k\}$ is a Cauchy sequence in~$\C^n$ converging to a point~$\sigma\in\de D$ by Lemma~\ref{dueotto}, and the assertion follows from Lemma~\ref{duetre}.
\end{proof}

The following general result proves Theorem~\ref{OJMz:}.(ii) in the strongly elliptic case.

\begin{lemma}\label{duedieci}
Let $D\Subset\C^n$ be a bounded strongly convex domain with $C^2$ boundary.
Let $f\in\Hol(D,D)$ be strongly elliptic with Wolff point $p\in D$. 
If $\sigma\in\de D$ is a boundary fixed point then 
$\beta_\sigma>1$.
\end{lemma}

\begin{proof} Since $p$ is a fixed point of $f$, we already know that
$$
\mlog\beta_\sigma=\liminf_{z\to\sigma}\bigl[k_D(z,p)-k_D\bigl(f(z),p\bigr)\bigr]\ge 0\;.
$$
Assume, by contradiction, that $\beta_\sigma=1$. Then Proposition~\ref{Julia} yields $f\bigl(E_p(\sigma,R)\bigr)\subseteq E_p(\sigma,R)$ for any $R>0$ because $\sigma$ is a boundary fixed point. 
Choose $R<1$ so that $p\notin\overline{E_p(\sigma,R)}$, and let $w\in\overline{E_p(\sigma,R)}$ be a point closest to~$p$ with respect to the Kobayashi distance. Since $f(w)\in\overline{E_p(\sigma,R)}$ this means that $k_D\bigl(f(w),p\bigr)\ge k_D(w,p)$, which is impossible because $w\ne p$ and $f$ is strongly elliptic.
\end{proof}

We conclude by proving Theorem~\ref{OJMz:}.(iii) in the strongly elliptic case.

\begin{lemma}
\label{quattro}
Let $D\Subset\C^n$ be a bounded strongly convex domain with $C^2$ boundary. Let $f\in\Hol(D,D)$ be strongly elliptic, with Wolff point $p\in D$. 
Let $\{z_k\}\subset D$ be a backward orbit for~$f$ with bounded Kobayashi step converging to~$\sigma\in\de D$. Then for every $q\in D$ 
there exists $M>0$ such that $z_k\in K_{q}(\sigma,M)$ eventually.
\end{lemma}

\begin{proof}
It suffices again to prove that there exists $M>1$ such that $\{z_k\}$ converges to $\sigma$ inside an admissible approach region $A(\sigma,M)$.

Without loss of generality, we can assume that $z_0\ne p$. We consider again $s_k>0$ defined by setting $-\mlog s_k=k_D(z_k,p)$. 
Thanks to \eqref{equsei}, there is a constant $0<c<1$ such that
\begin{equation}
s_{k+m} \le c^m s_k
\label{eq:alls}
\end{equation}
for all $k,m\ge 0$.

Now, \cite[Theorem 2.3.51, Theorem 2.3.52]{Aba1} yield constants $\widetilde C_1, \widetilde C_2 >0$ such that 
\begin{equation}
\widetilde C_1 d(z_j, \de D) \le s_j \le \widetilde C_2 d(z_j, \de D)
\label{eq:all}
\end{equation}
for all $j\in\N$, and so plugging this in \eqref{eqnove} we have
$$
\|z_k-z_{k+m}\| 
\le\frac{C}{1-\hat a}\frac{1}{1-c}\sqrt{s_k}
\le \frac{C}{1-\hat a}\frac{1}{1-c}\sqrt{\widetilde C_2} \sqrt{d(z_k, \de D)}
$$
for any $m\ge 0$ and $k$ large enough. Letting $m$ tend to infinity we then obtain
\begin{equation}
\|z_k-\sigma\|
\le M_1 \sqrt{d(z_k, \de D)}, 
\end{equation}
for some $M_1>1$.

On the other hand, up to translating the domain, without loss of generality we can assume that $D$ contains the origin. In particular, since $D$ is bounded and strongly convex we can replace $n_\sigma$ by $\sigma$ in the definition of $A(\sigma,M)$. Therefore, it suffices to prove that there exists $M_2>1$ such that 
$$
|\langle z_k-\sigma,\sigma\rangle| \le M_2 d(z_k, \de D)
$$
for $k$ large enough. But this follows by arguing as in the proof of Lemma \ref{tre} using $s_k$ instead of $t_k$, thanks to \eqref{eq:alls} and \eqref{eq:all}. Then taking $M= \max \{M_1, M_2\}$ we conclude the proof.
%
\end{proof}

This concludes the proof of Theorem \ref{OJMz:} in both cases.
\end{proof}

\bibliographystyle{alpha}

\end{document}